\documentclass[11pt]{article}
\usepackage{amssymb,latexsym,a4wide}
\newtheorem{corollary}{Corollary}
\newtheorem{remark}{Remark}
\newtheorem{lemma}{Lemma}
\newtheorem{proposition}{Proposition}
\newtheorem{theorem}{Theorem}

\newcommand{\ord}{ord}
\def\QQ{{\mathbb Q}}
\def\CC{{\mathbb C}}
\def\bull{\vrule height .9ex width .8ex depth -.1ex }

\begin{document}
\title{}
\author{D.~Grigoriev \\CNRS, IRMAR, Universit\'e de Rennes, Beaulieu,
35042, Rennes, France,\\ e-mail: dmitry.grigoryev@univ-rennes1.fr,\\
website: http://perso.univ-rennes1.fr/dmitry.grigoryev \and
F.~Schwarz
\\FhG, Institut SCAI, 53754 Sankt Augustin, Germany,\\
e-mail: fritz.schwarz@scai.fraunhofer.de\\
website: www.scai.fraunhofer.de/schwarz.0.html}
\date{}

\title{Non-holonomic overideals of differential ideals in two variables and absolute factoring}
\maketitle

\begin{abstract}
We study {\it non-holonomic} overideals of a (left) differential ideal
$J\subset F[d_x,\ d_y]$ in two variables
where $F$ is a differentially closed field of characteristic zero. The main result
states that a principal ideal $J=\langle P\rangle$ generated by an operator $P$ with
a separable {\it symbol} $symb(P)$ (being a homogeneous polynomial in two variables)
has a finite number of maximal non-holonomic overideals. This statement is extended
to non-holonomic ideals $J$ with a separable symbol. As an application we show that
in case of a second-order operator $P$ the ideal $\langle P\rangle$ has an infinite
number of maximal non-holonomic overideals iff $P$ is essentially ordinary. In case
of a third-order operator $P$  we give few sufficient conditions on
$\langle P\rangle$ to have a finite number of maximal non-holonomic overideals.

\end{abstract}

AMS Subject Classifications: 35A25, 35C05, 35G05

Keywords: differential non-holonomic overideals, Newton polygon, formal series
solutions

\section{Finiteness of a number of maximal non-holonomic overideals of an ideal
with a separable symbol}

Let $F$ be a differentially closed field (or universal in terms of
\cite{K}, \cite{S}) with derivatives $d_x,\ d_y$ and a linear
partial differential operator $P=\sum_{i,j} p_{i,j}d_x^id_y^j\in
F[d_x,\ d_y]$ be of an order $n$ (considering, e.~g. the field of
rational functions $\CC(x,y)$ as $F$ is a quite different issue).
The {\it symbol} $symb(P)= \sum_{i+j=n}p_{i,j}v^iw^j$ we treat as a
homogeneous polynomial in two variables of degree $n$. We call a
(left) ideal $I\subset F[d_x,\ d_y]$ {\it non-holonomic} if the
degree of its Hilbert-Kolchin polynomial $ez+e_0$ (in other words,
its {\it differential type} \cite{K}) equals 1. We study maximal
non-holonomic overideals of the principal ideal $ \langle P\rangle
\subset F[d_x,\ d_y]$ (obviously there is an infinite number of
maximal {\it holonomic} overideals of $\langle P\rangle$: for any
solution $u\in F$ of $Pu=0$ we get a holonomic overideal $\langle
d_x-u_x/u,d_y-u_y/u\rangle \supset \langle P\rangle$). We assume
w.l.o.g. that $symb(P)$ is not divisible by $d_y$ (otherwise one can
make a suitable transformation of the type $d_x\to d_x,\ d_y\to
d_y+bd_x,\ b\in F$, in fact choosing  $b$ from the subfield of
constants of $F$ would suffice).

Clearly, factoring an operator $P$ can be viewed as finding
principal overideals of  $\langle P\rangle$ and we refer to
factoring over a universal field $F$ as {\it absolute factoring}. We
mention also that overideals of an ideal in connection with Loewy
and primary decompositions were considered in \cite{GS07}.

Following \cite{G} consider a homogeneous polynomial ideal $symb(I)\subset F[v,w]$ and attach
a homogeneous polynomial $g=GCD(symb(I))$ to $I$. Lemma 4.1 \cite{G} states that $deg(g)=e$
(called also the {\it typical differential dimension} of $I$ \cite{K}). As above one can assume w.l.o.g.
that $w$ does not divide $g$.

We recall (see \cite{G05}, \cite{G}) that (Ore \cite{Bjork}) ring
$R=(F[d_y])^{-1}\ F[d_x,\ d_y]$ consists of fractions of the form
$\beta^{-1}r$ where $\beta\in  F[d_y],\ r\in  F[d_x,\ d_y]$. We also
recall that one can represent $R=F[d_x,\ d_y]\ (F[d_y])^{-1}$ and
two fractions are equal $\beta^{-1}r=r_1\beta^{-1}_1$ iff $\beta
r_1=r\beta_1$ \cite{G05}, \cite{G}.

For a non-holonomic ideal $I$ denote ideal $\overline{I}=RI\subset R$.
Since ring $R$ is left-euclidean (as well as right-euclidean) with respect to $d_x$ over
skew-field $(F[d_y])^{-1}\ F[d_y]$, we conclude
that ideal $\overline{I}$ is
principal, let $\overline{I}=\langle r\rangle$ for suitable $r\in F[d_x,\ d_y]
\subset R$ (cf. \cite{G}). Lemma 4.3 \cite{G} implies that $symb(r)=w^mg$ for a certain
integer $m\geq 0$ where $g$ is not divisible by $w$.

Now we expose a construction introduced in \cite{G}. For a family of elements $f_1,\dots,f_k\in F$
and rationals $1>s_2>\cdots>s_k>0$ we consider a $D$-module being a vector space over $F$ with a
basis $\{G^{(s)}\}_{s\in \QQ}$ where the derivatives of
$G^{(s)}=G^{(s)}(f_1,\dots,f_k;\ s_2,\dots,s_k)$ are defined as
$$d_{x_i}G^{(s)}=(d_{x_i}f_1)G^{(s+1)}+(d_{x_i}f_2)G^{(s+s_2)}+\cdots+(d_{x_i}f_k)G^{(s+s_k)}$$
\noindent for $i=1,2$ using the notations $d_{x_1}=d_x,d_{x_2}=d_y$.

Next we introduce series of the form
\begin{eqnarray}\label{0}
\sum_{0\leq i<\infty}h_iG^{(s-{i\over q})}
\end{eqnarray}
\noindent where $q$ is the least common multiple of the denominators of $s_2,\dots,s_k$ (one
can view (\ref{0}) as an analogue of Newton-Puiseux series for {\it non-holonomic} $D$-modules).

Theorem 2.5 \cite{G} states that for any (linear) divisor $v+aw$ of
$symb(P)$ and any $f_1\in F$ such that $(d_x+ad_y)f_1=0$ there
exists a solution of $P=0$ of the form (\ref{0}) (and conversely, if
(\ref{0}) is a solution of $P=0$ then $(d_x+ad_y)f_1=0$ for an
appropriate divisor $v+aw$ of $symb(P)$). Furthermore,
\mbox{Proposition 4.4 \cite{G}} implies that any solution of the
form (\ref{0}) of $r=0$ such that $(d_x+ad_y)f_1=0$ for suitable
$a\in F$ (or equivalently $d_yf_1\neq 0$) is also a solution of ideal
$I$ (then the appropriate linear form $v+aw$ is a divisor of $g$),
the inverse holds as well.

In \cite{GS04} we have designed an algorithm for factoring an operator $P$ is case when $symb(P)$
is separable. In particular, in this case there is only a finite number (less than $2^n$) of
different factorizations of $P$. Now we show a more general statement for overideals of $\langle
P \rangle$.

\begin{theorem}\label{finite}
Let $symb(P)$ be separable. Then there exists at most $n=ord(P)$
maximal non-holonomic overideals of $\langle P \rangle \subset
F[d_x,\ d_y]$. Moreover, if there exists a non-holonomic overideal
$I\supset \langle P\rangle$ with the attached polynomial
$g=GCD(symb(I))$ then there exists a unique non-holonomic
overideal maximal among ones with the attached polynomial equal $g$.
\end{theorem}

{\bf Proof}. Let a non-holonomic ideal $I\supset \langle
P \rangle$. Then $\beta P=r_1r$ for suitable $\beta\in F[d_y],\ r_1\in F[d_x,\ d_y]$ and a
polynomial $g=GCD(symb(I))$ attached to $I$ is a divisor of $symb(P)$.

We claim that for every pair of non-holonomic ideals $I_1,\
I_2\supset \langle P \rangle$ to which a fixed polynomial $g$ is
attached, to their sum $I_1+I_2$ also $g$ is attached. Indeed,
any solution of the form (\ref{0}) of $P=0$ such that $(v+aw)|g$, is
a solution of $r=0$ as well due to Lemma 4.2 \cite{G} (cf.
Proposition 4.4 \cite{G}) taking into account that $symb(P)$ is
separable, hence it is also a solution of $I$ as it was shown above
and by the same token is a solution of both $I_1$ and $I_2$ (in
particular $I_1+I_2$ is also non-holonomic). The claim is
established.

Thus among non-holonomic overideals  $I\supset \langle
P \rangle$ to which a given polynomial $g|symb(P)$ is attached, there is a unique maximal one.
Now take two maximal non-holonomic overideals $I,\ I'\supset \langle
P \rangle$ to which polynomials $g,\ g'$ are attached, respectively. Then $g,\ g'$ are
reciprocately prime. Indeed, if $v+aw$ divides both $g,\ g'$ then arguing as above one can verify
that (\ref{0}) is a solution of $I+I'$, i.~e. the latter ideal is non-holonomic which contradicts
to maximality of $I,\ I'$. Theorem is proved. \bull

\begin{corollary}\label{intersection}
Let $symb(P)$ be separable. Suppose that there exist maximal
non-holonomic overideals $I_1,\dots,I_l \supset \langle P\rangle$
such that for the respective attached polynomials $g_1,\dots,g_l$
the sum of their degrees $deg(g_1)+\cdots+deg(g_l)\geq n$. Then
$\langle P\rangle =I_1\cap \cdots \cap I_l$.
\end{corollary}

{\bf Proof}. As it was shown in the proof of Theorem~\ref{finite},
polynomials $g_j|symb(P),\ 1\leq j\leq l$ are pairwise reciprocately
prime, hence $g_1\cdots g_l =symb(P)$. Moreover it was established
in the proof of Theorem~\ref{finite} that every solution of $P=0$ of
the form (\ref{0}) such that $(d_x+ad_y)f_1=0$, is a solution of (a unique) $I_j$
for which $(u+aw)|g_j$, thus every solution of $P=0$ of the form (\ref{0})
is also a solution of $I_1\cap \cdots \cap I_l$. Therefore the typical differential dimension
of ideal $I_1\cap \cdots \cap I_l$ equals $n$ (cf. Lemma 4.1 \cite{G}).
On the other hand, any overideal of a principal ideal $\langle P\rangle$
of the same typical differential dimension coincides with $\langle P\rangle$
(one can verify it by comparing their Janet bases). \bull

\begin{remark}
One can extend Theorem~\ref{finite} to non-holonomic ideals $J$ such that homogeneous polynomial
$GCD(symb(J))$ is separable: namely, there exists a finite number of maximal non-holonomic
overideals $I\supset J$.
\end{remark}

\section{Non-holonomic overideals of a second-order linear partial differential operator}

In this Section we study the structure of overideals of $\langle P
\rangle$ when $n=ord(P)=2$. The case of separable $symb(P)$ is
covered by Theorem~\ref{finite}. Let $symb(P)$ be non-separable.
Then applying a transformation of the type $d_x\to b_1d_x+b_2d_y,\
d_y\to b_3d_x+b_4d_y$ for suitable $b_1,b_2,b_3,b_4\in F$ one can
assume w.l.o.g. that $P=d_y^2+p_1d_x+p_2d_y+p_3$ (it would be
interesting to find out when one can carry out these transformations
algorithmically). First let $p_1=0$. Then $P$ is essentially
ordinary (i.~e. becomes ordinary after a transformation as above)
and for any solution $u\in F$ of the equation $P=0$ we get a
non-holomonic overideal $\langle d_y-u_y/u\rangle \supset \langle P
\rangle$.

Now suppose that $p_1\neq 0$. Then $P$ is irreducible (see e.~g.
Corollary 7.1 \cite{G}). Moreover we claim that $\langle P\rangle$
has at most one maximal non-holonomic overideal. Let $I\supset
\langle P\rangle$ be a non-holonomic overideal. Choosing arbitrary
non-zero elements $b_1,b_2\in F$ denote the derivation
$d=b_1d_x+b_2d_y$. Similar to the proof of Theorem~\ref{finite}
there exists $r\in F[d,d_y]=F[d_x,d_y]$ such that $\langle r \rangle
=IR_1\subset R_1= (F[d])^{-1}\ F[d,d_y]$. Then $\beta P=r_1r$ for
suitable $\beta \in F[d],\ r_1\in F[d,d_y]$ and
$symb(r)=(b_1v+b_2w)^mg$ for an integer $m$ and $g|w^2$. If $g=1$
then $I$ cannot be non-holonomic because of Proposition 4.4 \cite{G}
(cf. above). If $g=w^2$ then similar to the proof of Corollary~\ref{intersection}
one can show that the only
non-holonomic overideal of $\langle P\rangle$ among ones to which
polynomial $w^2$ is attached, is just $\langle P\rangle$ itself.

It remains to consider the case $g=w$. Applying the Newton polygon
construction from \cite{G} to equation $r=0$ and a divisor $w$ of
$symb(r)$, one obtains a solution of the form (\ref{0}) of $r=0$
with $G=G(x)$, thereby it is a solution of $P=0$. On the other hand,
applying the Newton polygon construction from \cite{G} to equation
$P=0$, one gets at its first step $f_1=x$ and at the second step
$f_2$ which fulfils equation $(d_yf_2)^2+p_1=0$ and $f_2$
corresponds to the edge of the Newton polygon with endpoints
$(0,2),\ (1,0)$, so with the slope $1/2$. This provides a solution
of equation $P=0$ of the form (\ref{0}) with $G=G(x,f_2;\ 1/2)$,
therefore equation $P=0$ has no solutions of the form (\ref{0}) with
$G=G(x)$. The achieved contradiction shows that there are no
non-holonomic overideals $I$ with attached polynomial $w$, this completes the proof of
the claim.

Summarizing we can formulate (cf. \cite{GS04} for factoring $P$ over
not necessarily differentially closed fields)

\begin{proposition}\label{second}
Principal ideal $\langle P\rangle$ for a second-order operator
$P=d_y^2+p_1d_x+p_2d_y+p_3$ with non-separable $symb(P)$ has

i) no proper non-holonomic overideals in case $p_1\neq 0$;

ii) an infinite number of maximal non-holonomic overideals in case
$p_1=0$.
\end{proposition}

\section{On non-holonomic overideals of a third-order operator}

Now we study overideals of $\langle P\rangle$ where the order
$n=ord(P)=3$ (we mention that in \cite{G} an algorithm is designed
for factoring $P$). Due to Theorem~\ref{finite} it remains to
consider non-separable  $symb(P)$. We mention that few explicit
calculations for factoring $P$ are provided in \cite{GS08}.

\subsection{Symbol with two different linear divisors}

First let $symb(P)$ have two (linear) divisors and therefore one can assume w.l.o.g.
(see above) that
$w$ is its divisor of multiplicity 2 and $v$ is its divisor of multiplicity 1. One can
write
$$P=d_y^2d_x+p_0d_x^2+p_1d_y^2+p_2d_xd_y+p_3d_y+p_4d_x+p_5.$$

Suppose that $p_0\neq 0$. The Newton polygon construction from \cite{G} applied
to equation $P=0$ and to divisor
$w$ of $symb(P)$, yields a solution of the form (\ref{0}) of $P=0$ with $f_1=x$ at its
first step. At its second step the construction yields $f_2$ which fulfils equation
$(d_yf_2)^2+p_0=0$ and which corresponds to the edge of the Newton polygon with
endpoints $(1,2),\ (2,0)$, so with the slope $1/2$. This provides $G=G(x,f_2;\ 1/2)$ in
(\ref{0}).

Let a non-holonomic ideal $I\supset \langle P\rangle$. Choose $d=b_1d_x+b_2d_y$ for non-zero
$b_1,b_2\in F$. As in the previous Section there exists $r\in F[d,d_y]$ such that $\langle
r \rangle =R_1I\subset R_1=(F[d])^{-1}\ F[d,d_y]$. Then $\beta P=r_1r$ for suitable
$\beta\in F[d],\ r_1\in F[d,d_y]$. Rewrite $symb(r)=(b_1v+b_2w)^m g$ where $g|(vw^2)$. If
either $g=w^2$ or $g=v$, respectively, one can argue as in the proof of Theorem~\ref{finite}
and deduce that there can exist at most one maximal non-holonomic overideal of $\langle P\rangle$
with the property that attached to the overideal polynomial is either $w^2$ or $v$,
respectively. Similar to the proof of Corollary~\ref{intersection} one can verify that if
there exist maximal (non-holonomic) overideals $I_2,I_1\supset \langle P\rangle$ with attached
polynomials $w^2$ and $v$, respectively, then $\langle P\rangle = I_1\cap I_2$. As in
Theorem~\ref{finite} the existence of a maximal overideal with the attached polynomial $w^2$
(or respectively, $v$) follows from the existence of any non-holonomic overideal with the
attached polynomial $w^2$ (or respectively, $v$).

If either $g=w$ or $g=vw$ then applying the Newton polygon construction from \cite{G} to
equation $r=0$ and divisor $w$ of $symb(r)$, one obtains a solution of $r=0$ (and thereby,
of $P=0$ due to Lemma 4.2 \cite{G}) of the form (\ref{0}) with $G=G(x)$ which contradicts
to the supposition $p_0\neq 0$ (see above).
Thus, in case  $p_0\neq 0$ ideal  $\langle P\rangle$ has a finite (less or equal than 2)
number of maximal
non-holonomic overideals (similar to Theorem~\ref{finite}).

When $p_0=0$ this is not always true, say for $P=(d_x+b)(d_y^2+b_3d_y+b_4)$
(cf. case $n=2$ in the previous Section).
It would be interesting to clarify for which $P$ this is still true.

\subsection{Symbol with a unique linear divisor}

Now we consider the last case when $symb(P)$ has a unique linear divisor with multiplicity
3. As above one can assume w.l.o.g. that $symb(P)=w^3$, so
$$P=d_y^3+p_0d_x^2+p_1d_y^2+
p_2d_xd_y+p_3d_y+p_4d_x+p_5.$$
\noindent Keeping the notations we get $\langle
r \rangle=R_1I$ and $\beta P=r_1r$. Then $symb(r)=(b_1v+b_2w)^m g$ where $g|w^3$. If
$g=w^3$ then arguing as in the proof of Corollary~\ref{intersection} we deduce that the only
non-holonomic overideal  of $\langle P \rangle$ to which polynomial $w^3$ is
attached, is just $\langle P \rangle$ itself.

Let $g|w^2$. Applying the Newton polygon construction from \cite{G} to equation $r=0$ and
linear divisor $w$ of $symb(r)$ one gets a solution of $r=0$ (and thereby of $P=0$)
with either
$G=G(x)$ or $G=G(x,f_2;\ 1/2)$ where $d_yf_2\neq 0$ (cf. above).

Application of the Newton polygon construction from \cite{G} to equation $P=0$ (and unique
linear divisor $w$ of $symb(P)$) at its first step provides $f_1=x$. The second step
requires a trial of cases. First let $p_0\neq 0$. Then the second step yields $f_2$ which
fulfils equation $(d_yf_2)^3+p_0=0$ and which corresponds to the edge of the Newton polygon
with endpoints $(0,3),\ (2.0)$, so with the slope $2/3$. Thus we obtain a solution of the
form (\ref{0}) with $G=G(x,f_2,\dots;\ 2/3,\dots)$, hence $\langle P\rangle$ in case $p_0\neq 0$ has
no non-holonomic overideals with attached polynomial $g$ being a divisor of $w^2$ (see
above).

Now assume that $p_0=0$ and $p_2\neq 0$. Then the second step provides solutions of
$P=0$ of the form (\ref{0}) with two different possibilities. Either the
Newton polygon construction chooses the vertical edge with endpoints $(1,1),\ (1,0)$
as a leading edge at the second step, then it terminates at the second step yielding a
solution  of the form (\ref{0}) with
$G=G(x)$ (we recall that in the construction from Section 2 \cite{G} only edges with
non-negative
slopes are taken as leading ones and the construction terminates while taking a vertical
edge, so with the slope $0$, as a leading one, in particular the edge with
endpoints
$(1,1),\ (1,0)$ is taken as a leading one regardless of whether the coefficient at point
$(1,0)$ vanishes). As the second possibility the construction yields  a solution  of the
form (\ref{0}) with $G=G(x,f_2,\dots;\ 1/2,\dots)$ where $f_2\neq 0$ fulfils equation
$(d_yf_2)^3+p_2d_yf_2=0$ corresponding to the edge of the Newton polygon with endpoints
$(0,3),\ (1,1)$, so with the slope $1/2$. One can suppose w.l.o.g. that the
 Newton polygon construction terminates at its third step (thereby $G=G(x,f_2;\ 1/2)$),
otherwise $\langle P\rangle$ cannot have a non-holonomic overideal to which a divisor $g$
of $w^2$ is attached (see above).

If $g=w^2$ then any solution $H_2$ of $P=0$ of the form (\ref{0}) with $G=G(x,f_2;\ 1/2)$
is a solution of $r=0$ because otherwise $rH_2\neq 0$, being also of the form (\ref{0})
with $G=G(x,f_2;\ 1/2)$,
cannot be a solution of $r_1=0$ taking into account that $symb(r_1)$ does not divide on $w^2$
 (cf. Lemma 4.2 \cite{G}). Else if $g=w$ then $rH_2\neq 0$ (again taking into account that
$symb(r)$ does not divide on $w^2$) and therefore $r_1(rH_2)=0$. Hence for a solution
$H_1$ of $P=0$ of the form (\ref{0}) with $G=G(x)$ (see above) we have $rH_1=0$ since
otherwise $rH_1$ being also of the form (\ref{0}) with $G=G(x)$ cannot be a solution of
$r_1=0$ (again cf. Lemma 4.2 \cite{G}). Then arguing as in the proof of Theorem~\ref{finite}
one concludes that in case $p_0=0$ and $p_2\neq 0$ ideal $\langle P\rangle$ can have at
most two maximal non-holonomic overideals (with attached polynomials $w$ and $w^2$,
respectively).
Similar to the proof
of Corollary~\ref{intersection} (cf. the preceding Subsection) one can verify that if there
exist maximal (non-holonomic) overideals $I_1,I_2\supset \langle P\rangle$ with attached
polynomials $w$ and $w^2$, respectively, then $\langle P\rangle = I_1\cap I_2$. As in
Theorem~\ref{finite} the existence of a maximal overideal with the attached polynomial $w$ (or
respectively, $w^2$) follows from the existence of any non-holonomic overideal with the
attached polynomial $w$ (or respectively, $w^2$).

Furthermore, let $p_0=p_2=0,\ p_4\neq 0$. Then as in case $p_0\neq 0$ we argue that the
second step of the Newton polygon construction applied to equation $P=0$ yields $f_2$ which
fulfils equation
$(d_yf_2)^3+p_4=0$ and which corresponds to the leading edge of the Newton polygon with
endpoints $(0,3),\ (1,0)$, so with the slope $1/3$. Thus the Newton polygon construction
yields a solution of $P=0$ of the form (\ref{0}) with $G=G(x,f_2,\dots;\ 1/3,\dots)$ and
again
$\langle P\rangle$ in case $p_0=p_2=0,\ p_4\neq 0$ under consideration has no non-holonomic
overideals with an attached polynomial being a divisor of $w^2$.

Finally, when $p_0=p_2=p_4=0$ ideal $\langle P=d_y^3+p_1d_y^2+p_3d_y+p_5\rangle$
has an infinite number of maximal non-holonomic overideals (similar to the second-order
case $P=d_y^2+p_3d_y+p_5$, see above). Summarizing we conclude with the following

\begin{proposition}
Let $P$ be a third-order operator with a non-separable $symb(P)$.

i) When $symb(P)$ has two different (linear) divisors (one of which of multiplicity 2) then
we can assume w.l.o.g. that $P=d_y^2d_x+p_0d_x^2+p_1d_y^2+p_2d_xd_y+p_3d_y+p_4d_x+p_5$.
If $p_0\neq 0$ then $\langle P\rangle$ has at most two maximal non-holonomic overideals.
Moreover if there exist two different maximal non-holonomic overideals $I_1,I_2\supset
\langle P\rangle$ then $\langle P\rangle = I_1\cap I_2$;

ii) when $symb(P)$ has a single linear divisor (of multiplicity 3) then we can assume
w.l.o.g. that $P=d_y^3+p_0d_x^2+p_1d_y^2+
p_2d_xd_y+p_3d_y+p_4d_x+p_5$. If either $p_0\neq 0$, either $p_2\neq 0$ or $p_4\neq 0$
then $\langle P\rangle$ has at most two maximal non-holonomic overideals.
Moreover if there exist two different maximal non-holonomic overideals $I_1,I_2\supset
\langle P\rangle$ then $\langle P\rangle = I_1\cap I_2$.
Otherwise $\langle P=d_y^3+p_1d_y^2+p_3d_y+p_5\rangle$ has an infinite number of
maximal non-holonomic overideals.
\end{proposition}

It is a challenge to design an algorithm which produces non-holonomic overideals of a
given differential ideal $J\subset F[d_x,d_y]$.

\section*{Appendix. Explicit formulas for Laplace transformation}

We exhibit a short exposition and explicit formulas for the Laplace
transformation \cite{Goursat}.

Let $Q=d_{xy}+ad_x+bd_y+c$ be a second-order operator and
$L_n=\sum_{0\leq i\le n} l_id_x^i$ its "Laplace divisor" of order
$n$, in particular $Q,\ L_n$ form a Janet basis, hence
\begin{eqnarray}\label{1.1}
PQ=(d_y+a)L_n
\end{eqnarray}
\noindent for a suitable $P=\sum_{0\le i\le n-1} p_id_x^i$ (the
latter shape of $P$ we obtain comparing the highest terms which
divide on $d_x^n$ in (\ref{1.1})). Comparing the highest terms in
(\ref{1.1}) which divide on $d_y$, we get that $L_n=P(d_x+b)$. Thus
\begin{eqnarray}\label{1.2}
PQ=(d_y+a)P(d_x+b).
\end{eqnarray}
\noindent We have $Q\neq (d_y+a)(d_x+b)$ iff $0\neq ab+b_y-c=:K_0$.
\begin{lemma}\label{transformation}
If $K_0\neq 0$ then there are unique $B,\ C$ such that
\begin{eqnarray}\label{1.3}
(d_x+B)Q=(d_{xy}+ad_x+Bd_y+C)(d_x+b)
\end{eqnarray}
\end{lemma}
{\bf Proof}. (\ref{1.3}) is equivalent to an algebraic linear in
$B,\ C$ system
\begin{eqnarray}\label{1.4}
aB-C=b_y+ab-a_x-c,\quad (c-b_y)B-bC=b_{xy}+ab_x-c_x
\end{eqnarray}
\bull

Therefore (\ref{1.2}) holds iff $P=P_1(d_x+B)$ (by means of dividing
$P$ by $d_x+B$ with remainder). Substituting the latter equality to
(\ref{1.2}) and making use of (\ref{1.3}) we obtain the equality
\begin{eqnarray}\label{1.5}
P_1(d_{xy}+ad_x+Bd_y+C)=(d_y+a)P_1(d_x+B).
\end{eqnarray}
\noindent Now (\ref{1.5}) is similar to (\ref{1.2}) but with the
order $\ord(P_1)=\ord(P)-1=n-1$ and a new second-order operator
$Q_1=d_{xy}+ad_x+Bd_y+C$. Continuing this way we get the Laplace
transformation with $K_1=aB+B_y-C$ etc.

More uniformly denote $b_0:=b,\ c_0:=c$, then $b_1:=B,\ c_1:=C,\
b_2,\ c_2$ etc. obtained from Lemma~\ref{transformation}. Denote
$$K_i:=ab_i+(b_i)_y-c_i,\ Q_i:=d_{xy}+ad_x+b_id_y+c_i.$$
\begin{corollary}
There exists $L_n$ satisfying (\ref{1.1}) iff for the minimal $m$
such that $K_m=0$ we have $m\le n$. In this case
\begin{eqnarray}\label{1.6}
L_n=P_{n-m}(d_x+b_{m-1})\cdots (d_x+b_0)
\end{eqnarray}
\noindent where $P_{n-m}=\sum_{0\le i\le n-m} p_id_x^i$ is an
arbitrary operator of the order $n-m$ which fulfils
\begin{eqnarray}\label{1.7}
P_{n-m}(d_y+a)=(d_y+a)P_{n-m}.
\end{eqnarray}
\noindent For any order $n-m\ge 0$ such an operator $P_{n-m}$
exists. The pair $Q,\ L_n$ constitutes a Janet basis of the ideal
$\langle Q,\ L_n\rangle$. The ideal $\langle Q,\ L_m\rangle$ is the
(unique) maximal non-holonomic overideal of $\langle Q \rangle$
which corresponds to a divisor $y$ of $symb(Q)=xy$ (see
Theorem~\ref{finite}).
\end{corollary}
{\bf Proof}. Applying Laplace transformations as above, if $m>n$ we
don't get a solution of (\ref{1.1}) after $n$ steps since
(\ref{1.2}) with $PQ_n=(d_y+a)P(d_x+b_n)$ would not have a solution
with $P$ of the order $0$. If $m\le n$ then successively following
Laplace transformations we arrive to (\ref{1.6}) in which
(\ref{1.7}) is obtained from equality $PQ_m=(d_y+a)P(d_x+b_m)$ (see
(\ref{1.2})) and taking into account that $K_m=0$. \bull

 \vspace{4mm} {\bf Acknowledgement}. The first author is
grateful to the Max-Planck Institut f\"ur Mathematik, Bonn for its
hospitality while writing this paper.


\begin{thebibliography}{xx}



\bibitem{Bjork}
J.~E.~Bj\"ork, Rings of differential operators, North-Holland, 1979.





\bibitem{Goursat}
E.~Goursat, Le\c con sur l'int\'egration des \'equations aux
d\'eriv\'ees partielles, vol.~I, II, A.~Hermann, 1898.





\bibitem{G05} D.~Grigoriev, Weak B\'{e}zout Inequality for D-Modules, {\em J. Complexity} {\bf 21} (2005), 532-542.

\bibitem{G} D.~Grigoriev, Newton-Puiseux series for non-holonomic D-modules
and factoring linear partial differential operators

\bibitem{GS04} D.~Grigoriev, F.~Schwarz, Factoring and solving linear
partial differential equations, {\em Computing} {\bf 73} (2004), 179-197.

\bibitem{GS07} D.~Grigoriev, F.~Schwarz, Loewy and primary decomposition of D-Modules,
{\em Adv. Appl. Math.} {\bf 38} (2007), 526-541.

\bibitem{GS08} D.~Grigoriev, F.~Schwarz, Loewy decomposition of
linear third-order PDE's in the plane, Proc. Intern. Symp. Symbolic,
Algebr. Comput., ACM Press, 2008.

\bibitem{K} E.~Kolchin, Differential Algebra and Algebraic Groups,
Academic Press, New York, 1973.

\bibitem{S} M.~van der Put, M.~Singer, Galois theory of linear
differential equations, Grundlehren der Mathematischen
Wissenschaften, {\bf 328}, Springer, 2003.





















\end{thebibliography}
\end{document}